\documentclass{proc-l}

\newtheorem{thm}{Theorem}

\theoremstyle{remark}

\theoremstyle{definition}

\title{Minimizing Euler characteristics of symplectic four-manifolds}
\author{D.~Kotschick}
\address{Mathematisches Institut, Ludwig-Maximilians-Universit\"at M\"unchen,
Theresienstr.~39, 80333 M\"unchen, Germany}
\email{dieter@member.ams.org}
\thanks{I am grateful to P.~Kirk for pointing out the question in~\cite{BK} that is answered here.}
\date{May 3, 2005}

\begin{document}

\begin{abstract}
We prove that the minimal Euler characteristic of a closed symplectic four-manifold with given fundamental group is often
much larger than the minimal Euler characteristic of almost complex closed four-manifolds with the same fundamental group.
In fact, the difference between the two is arbitrarily large for certain groups.
\end{abstract}

\maketitle

It was first proved by Dehn~\cite{Dehn} that every finitely presentable group $\Gamma$ can be realized as the fundamental 
group of a closed oriented smooth $4$-manifold. Taking the minimum over the Euler characteristics of
all such manifolds one obtains an interesting numerical invariant $q^{DIFF}(\Gamma)$ of finitely presentable groups, see
for example~\cite{Gromov,HW,tani}. As mentioned in~\cite{tani}, there are geometric variants $q^{GEO}(\Gamma)$ of this 
definition, obtained by minimizing the Euler characteristic only over those $4$-manifolds with fundamental group $\Gamma$
which carry a specified geometric structure. One trivially has 
$$
q^{DIFF}(\Gamma)\leq q^{GEO}(\Gamma)
$$
for all geometric structures. Moreover, the inequality is often strict.

For a simple example of a geometric invariant, consider almost complex $4$-manifolds. Every finitely presentable group is the fundamental 
group of an almost complex $4$-manifold~\cite{BLMS}, but the minimal Euler characteristic over almost complex $4$-manifolds is strictly 
larger than $q^{DIFF}(\Gamma)$ for many $\Gamma$. Nevertheless, in this case it is easy to see that the difference between the smooth and 
geometric invariants is universally bounded independently of $\Gamma$, compare~\cite{BLMS}.

The purpose of this paper is to show that in the symplectic category this boundedness fails. Recall that Gompf~\cite{Gompf}
proved that every finitely presentable $\Gamma$ can be realised as the fundamental group of a closed symplectic $4$-manifold. 
Thus we can define $q^{SYMP}(\Gamma)$ to be the minimal Euler characteristic of a closed symplectic $4$-manifold with fundamental 
group $\Gamma$. Then we have
\begin{thm}\label{t:main}
For every $c>0$ there exists a finitely presentable group $\Gamma$ satisfying
$$
q^{SYMP}(\Gamma)\geq q^{DIFF}(\Gamma) + c \ .
$$
\end{thm}
\begin{proof}
We shall use the sequence $F_r$ of free groups of rank $r$. It suffices to show that the difference
$$
q^{SYMP}(F_r) - q^{DIFF}(F_r) 
$$
grows linearly with the rank $r$. We know from~\cite{tani} that $q^{DIFF}(F_r) =-2(r-1)$, because, on the one hand, this 
value is the obvious lower bound $2-2b_1$ for the Euler characteristic of any closed $4$-manifold with fundamental group
$F_r$, and, on the other hand, this value is realized by the connected sum of $r$ copies of $S^1\times S^3$.

To estimate $q^{SYMP}(F_r)$ let $X$ be a closed symplectic $4$-manifold with fundamental group $F_r$ and with minimal Euler characteristic.
The minimality of the Euler characteristic implies that $X$ is symplectically minimal in the sense that it contains no symplectically
embedded $(-1)$-spheres. Let us assume for the moment that the positive part $b_2^+(X)$ of the intersection form of $X$ is
strictly larger than $1$, then a result of Taubes~\cite{Taubes} implies $c_1^2(X)\geq 0$, see also~\cite{Bourbaki}. 
We expand this inequality as follows:
$$
0\leq c_1^2=2\chi+3\sigma=4-4b_1+5b_2^+-b_2^-\leq 4-4b_1+5b_2^+ \ .
$$
This yields $b_2^+\geq\frac{4}{5}(b_1-1)$, and thus
$$
\chi=2-2b_1+b_2\geq 2-2b_1+b_2^+\geq -\frac{6}{5}(b_1-1) \ .
$$
Therefore we have
\begin{equation}\label{bound}
q^{SYMP}(F_r)\geq -\frac{6}{5}(r-1) \ ,
\end{equation}
showing that the difference $q^{SYMP}(F_r) - q^{DIFF}(F_r)$ grows linearly with $r$.

It remains to remove the assumption $b_2^+(X)>1$. As $X$ is symplectic, the only other possibility is $b_2^+(X)=1$. If this happens, 
consider a $d$-fold covering $X_d$ of $X$, with $d>1$. This is symplectic with free fundamental group of rank $1+d(r-1)$. The 
multiplicativity of the signature and of the Euler characteristic in finite coverings then imply $b_2^+(X_d)=db_2^+(X)=d>1$. We can not 
apply Taubes's inequality to $X_d$ because a priori we do not know that $X_d$ is symplectically minimal. Instead of proving this, we 
argue as follows. The minimal model $Y_d$ of $X_d$ has the same $b_1$ and $b_2^+$ as $X_d$. Taubes's inequality $c_1^2\geq 0$ 
applied to $Y_d$ gives
$$
0\leq c_1^2(Y_d)\leq 4-4b_1(Y_d)+5b_2^+(Y_d)=d(9-4r) \ .
$$
It follows that $r\leq 2$. In the cases $r\leq 1$, inequality~\eqref{bound} is trivial. In the case $r=2$ it reduces to $q^{SYMP}(F_2)\geq -1$,
which is true because in this case $\chi (X)= 2-2b_1(X)+b_2(X) = -2+b_2(X)\geq -1$.
\end{proof}
This result was motivated by the recent paper of Baldridge and Kirk~\cite{BK}, concerned with a systematic study of $q^{SYMP}(\Gamma)$.
The lower bounds for $q^{SYMP}(\Gamma)$ given in~\cite{BK} are never better than $q^{DIFF}(\Gamma)+2$, because only the condition
$b_2^+\geq 1$ and the existence of almost complex structures on symplectic manifolds are used.

It turns out that the bound~\eqref{bound} holds in almost complete generality:
\begin{thm}
Let $\Gamma$ be a finitely presentable group. The inequality
\begin{equation}\label{bound2}
q^{SYMP}(\Gamma)\geq -\frac{6}{5}(b_1(\Gamma)-1) 
\end{equation}
holds for $\Gamma$ if and only if $\Gamma$ is not the fundamental group of a closed oriented surface of genus $\geq 2$.
\end{thm}
\begin{proof}
First of all, if $\Gamma$ is the fundamental group of a closed oriented surface of genus $g\geq 2$, then it was proved in~\cite{tani}
that $q^{DIFF}(\Gamma)=4(1-g)=2(2-b_1(\Gamma))$. The manifold $S^2\times\Sigma_g$ realizes the minimum and is symplectic,
so that $q^{SYMP}(\Gamma)=2(2-b_1(\Gamma))<\frac{6}{5}(1-b_1(\Gamma))$.

Suppose now that $\Gamma$ is not a surface group. If a symplectic manifold with fundamental group $\Gamma$ realizing the 
smallest possible Euler characteristic has $b_2^+>1$, then Taubes's inequality $c_1^2\geq 0$ for minimal symplectic manifolds 
with $b_2^+>1$ implies~\eqref{bound2}, as in the proof of Theorem~\ref{t:main}. If the symplectic minimizer for $\Gamma$ has
$b_2^+=1$, then for arbitrary $\Gamma$ we may not be able to use covering tricks as in the proof of Theorem~\ref{t:main}.
However, because $\Gamma$ is not a surface group, our manifold cannot be ruled. Therefore we can use Liu's extension~\cite{Liu}
of Taubes's inequality to minimal non-ruled symplectic manifolds with $b_2^+=1$ to reach the same conclusion as before.
\end{proof}

Gompf~\cite{Gompf} asked whether a non-ruled symplectic $4$-manifold necessarily has non-negative Euler characteristic. 
This question is still open. A positive answer would of course provide a vast generalization of the results proved here. If a 
finitely presentable group $\Gamma$ satisfies $q^{DIFF}(\Gamma)<0$, then one knows a lot of its properties. For example, 
$\Gamma$ cannot embed non-trivially in itself with finite index, it is non-amenable, and has a subgroup of finite index surjecting 
onto $F_2$, see~\cite{Gromov,tani}. Thus there are many group-theoretic constraints for a negative answer to Gompf's question.

\bigskip

\bibliographystyle{amsplain}

\bigskip

\end{document}